\documentclass[10pt,a4paper]{article}
\usepackage{graphics, graphicx}     

\usepackage{cite}
\usepackage{amsmath}
\usepackage{amssymb}
\usepackage{amsthm}
\theoremstyle{plain}
\newtheorem{Theorem}{Theorem}[section] %
\newtheorem{Lemma}{Lemma}[section]

\theoremstyle{definition}
\newtheorem{Remark}{Remark}[section]

\theoremstyle{definition}

\newtheorem{Problem}{Problem}[section]

{\par\noindent{\it Proof of}} 
{\hfill$\vspace{5mm}\scriptstyle\blacksquare$} 

\numberwithin{equation}{section} 
\numberwithin{figure}{section} 
\numberwithin{table}{section} 

\begin{document}

\setcounter{page}{1}

\markboth{M.I. Isaev, R.G. Novikov}{New global stability estimates for monochromatic inverse acoustic  scattering}

\title{New global stability estimates for monochromatic inverse acoustic scattering}
\date{}
\author{ { M.I. Isaev and R.G. Novikov}}

\maketitle
\begin{abstract}
We give new global stability estimates for monochromatic inverse acoustic scattering. 
These estimates essentially improve estimates of [P. H\"ahner, T. Hohage, 
SIAM J. Math. Anal., 33(3), 2001, 670-685] and can be considered as a solution of an open problem
formulated in the aforementioned work. 
\end{abstract}

\section{Introduction}
We consider the equation
\begin{equation}\label{eq} 
	\Delta \psi + \omega^2 n(x) \psi =0, \ \ \  x \in  \mathbb{R}^3, \ \omega>0,
\end{equation}
where
\begin{equation}\label{eq_c}
	\begin{aligned}
		 (1-n) \in \mathbb{W}^{m,1}(\mathbb{R}^3) \text{ for some } m > 3, \\ 
		 \mbox{Im}\, n(x) \geq 0, \ \ \ x\in \mathbb{R}^3,\\
		 \mbox{supp}\, (1 - n) \subset B_{r_1}  \text{ for some } r_1 > 0,
	\end{aligned}
\end{equation}
where 
 $\mathbb{W}^{m,1}(\mathbb{R}^3)$ denotes the standart Sobolev space on $\mathbb{R}^3$ (see formula \eqref{2.5} of Section 2 for details), 
 $B_r = \{x\in \mathbb{R}^3: |x|<r\}$. 
 
 We interpret \eqref{eq} as the stationary acoustic equation at frequency $\omega$ in an inhomogeneous 
 medium with refractive index $n$.
 
In addition, we consider the Green function $G^+(x,y,\omega)$ for the operator $\Delta+\omega^2 n(x)$ with 
  the Sommerfeld radiation condition:
\begin{equation}\label{G+def}
\begin{aligned}
	\left(\Delta + \omega^2 n(x)\right)G^+(x,y,\omega) = \delta(x-y),\\
	\lim\limits_{|x|\rightarrow \infty} |x|
	\left(
		\frac{\partial G^+}{\partial |x|}(x,y,\omega) - i\omega G^+(x,y,\omega)
	\right)=0,\\
	\text{uniformly for all directions $\hat{x} = x/|x|$},
	\\
	 x,y\in \mathbb{R}^3, \ \omega>0.
\end{aligned}
\end{equation} 
It is know that, under assumptions \eqref{eq_c}, the function $G^+$ is uniquely specified
by \eqref{G+def}, see, for example, \cite{HH2001}, \cite{CK1998}.

We consider, in particular, the following near-field inverse scattering problem for equation (\ref{eq}): 

\begin{Problem}
 Given $G^+$ on $\partial B_r\times \partial B_r$ for some fixed $\omega>0$ and $r>r_1$, find $n$ on $B_{r_1}$. 
\end{Problem}

We consider also the solutions  $\psi^{+}(x,k)$, $x\in \mathbb{R}^3, \	k\in \mathbb{R}^3, \ k^2 =\omega^2$, of
equation \eqref{eq} specified by the following asymptotic condition:
	\begin{equation}\label{Psi_scattering}
	\begin{array}{c} \displaystyle
		\psi^+(x,k) = e^{ikx} - 2\pi^2 \frac{e^{i|k||x|}}{|x|} f\left(k,|k|\frac{x}{|x|}\right)+ o\left(\frac{1}{|x|}\right) \\\\
		\displaystyle
		\text{as } \ |x| \rightarrow \infty \left( \text{uniformly in } \frac{x}{|x|}\right),
	\end{array}
\end{equation}	 
	with some a priory unknown $f$.

The function $f$ on 
$
	\mathcal{M}_\omega = \{k\in \mathbb{R}^3, l \in \mathbb{R}^3: k^2 =l^2 = \omega^2\}
$
arising in (\ref{Psi_scattering}) is  the classical scattering amplitude 
for equation (\ref{eq}). 

In addition to Problem 1.1, we consider also the following far-field inverse scattering  problem for equation (\ref{eq}):

\begin{Problem}
	Given $f$ on $\mathcal{M}_\omega$ for some fixed $\omega>0$, find $n$ on $B_{r_1}$. 
\end{Problem}

In \cite{Berezanskii1958} it was shown that the near-field data of Problem 1.1 are uniquely determined by the far-field data of Problem 1.2 and vice versa.

Global uniqueness for Problems 1.1 and 1.2 was proved for the first time in \cite{Novikov1988}; in addition, this proof is constructive. For more information
on reconstruction methods for Problems 1.1 and 1.2 see \cite{ABR2008}, \cite{HH2001}, \cite{Nachman1988}, \cite{Novikov1988}, \cite{Novikov2005+}, \cite{NS2012} 
and references therein.

Problems 1.1  and 1.2 can be also considered as  examples of ill-posed problems: 
see \cite{LR1986}, \cite{BK2012} for an introduction to this theory.

The main results of the present article consist of the following two theorems:
\begin{Theorem}\label{Theorem_1.1} Let $C_n>0$, $r>r_1$ be fixed constants.  
Then there exists a positive constant $C$ (depending only on $m$, $\omega$, $r_1$, $r$ and $C_n$) such that
for all refractive indices $n_1$, $n_2$ satysfying 
$\|1-n_1\|_{\mathbb{W}^{m,1}(\mathbb{R}^3)}, \|1-n_2\|_{\mathbb{W}^{m,1}(\mathbb{R}^3)}<C_n$,
{\rm $\mbox{supp}\, (1-n_1), \mbox{supp}\, (1-n_2)\subset B_{r_1}$}, 
the following estimate holds:
\begin{equation}\label{eq_t1}
	||n_1 - n_2||_{\mathbb{L}^\infty(\mathbb{R}^3)} \leq 
	C\left(\ln\left(3+\delta ^{-1}\right)\right)^{-s}, \ \ s = \frac{m-3}{3},
\end{equation}
where   $\delta  = ||G^+_1 - G^+_2||_{\mathbb{L}^{2}(\partial B_r \times \partial B_r)}$ and $G_1^+$, $G_2^+$ 
are the near-field scattering data for the refractive indices $n_1$, $n_2$, respectively, at fixed frequency $\omega$. 
\end{Theorem}

\begin{Remark}
We recall that if $n_1,n_2$ are refractive indices satisfying \eqref{eq_c}, then
 $
	G^+_1 - G^+_2 \text { is bounded  in } \mathbb{L}^{2}(\partial B_r \times \partial B_r) \text{ for any $r>r_1$,}
$
  where $G_1^+$ and $G_2^+$ 
are the near-field scattering data for the refractive indices $n_1$ and $n_2$, respectively, at fixed frequency $\omega$, 
see, for example, Lemma 2.1 of \cite{HH2001}.	
\end{Remark}

\begin{Theorem}\label{Theorem_1.2} Let $C_n>0$ and $0<\epsilon <\frac{m-3}{3}$ be fixed constants.  
Then there exists a positive constant $C$ (depending only on $m$, $\epsilon$, $\omega$, $r_1$ and $C_n$) such that
for all refractive indices $n_1$, $n_2$ satysfying 
$\|1-n_1\|_{\mathbb{W}^{m,1}(\mathbb{R}^3)}, \|1-n_2\|_{\mathbb{W}^{m,1}(\mathbb{R}^3)}<C_n$,
{\rm $\mbox{supp}\, (1-n_1), \mbox{supp}\, (1-n_2)\subset B_{r_1}$}, 
the following estimate holds:
\begin{equation}\label{eq_t2}
	||n_1 - n_2||_{\mathbb{L}^\infty(\mathbb{R}^3)} \leq 
	C\left(\ln\left(3+\delta ^{-1}\right)\right)^{-s+\epsilon}, \ \ s = \frac{m-3}{3},
\end{equation}
where   $\delta  = ||f_1 - f_2||_{\mathbb{L}^{2}(\mathcal{M}_\omega)}$ and $f_1$, $f_2$ 
denote the scattering amplitudes for the refractive indices $n_1$, $n_2$, respectively, at fixed frequency $\omega$. 
\end{Theorem}

For some regularity dependent $s$ but always smaller than $1$ the stability estimates of Theorems \ref{Theorem_1.1} and \ref{Theorem_1.2}
were proved in \cite{HH2001}. Possibility of estimates \eqref{eq_t1}, \eqref{eq_t2} with $s>1$ was formulated in \cite{HH2001} as 
an open problem, see page 685 of \cite{HH2001}. Our estimates \eqref{eq_t1}, \eqref{eq_t2} with $s= \frac{m-3}{3}$ give a solution of this problem. 
Apparently, using the methods of \cite{Novikov2009}, \cite{Novikov2011} estimates \eqref{eq_t1}, \eqref{eq_t2} can be proved for $s=m-3$.
For more information on stability estimates for Problems 1.1 and 1.2 see \cite{HH2001}, \cite{IsaevFunc}, \cite{Stefanov1990} and references therein.
In particular, as a corollary of \cite{IsaevFunc} estimates \eqref{eq_t1}, \eqref{eq_t2} can not be fulfilled, in general, for $s>\frac{5m}{3}$.

The proofs of Theorem \ref{Theorem_1.1} and \ref{Theorem_1.2} are given in Section 3. These proofs use, in particular: 
\begin{enumerate}
\item Properties of the
Faddeev functions for equation \eqref{eq} considered as the Schr\"odinger equation at fixed energy $E=\omega^2$, see Section 2. 
\item The results of \cite{HH2001} consisting in Lemma \ref{Lemma_3.1} and in reducing (via Lemma 3.2) estimates of the form \eqref{eq_t2} 
for Problem 1.2 to estimates of the form  \eqref{eq_t1} for Problem 1.1.
\end{enumerate}

In addition in the proofs of Theorem  \ref{Theorem_1.1} and \ref{Theorem_1.2} we combine some of the aforementioned ingredients in a similar way
with the proof of stability estimates of \cite{IN2012++}. 

\section{Faddeev functions}
We consider \eqref{eq} as the Schr\"odinger equation at fixed energy $E=\omega^2$: 
\begin{equation}\label{eq_2.1}
	-\Delta\psi + v(x)\psi = E \psi, \ \ \ x\in \mathbb{R}^3,
\end{equation}
where
$
	v = \omega^2(1-n)$,  $E = \omega^2$.

For equation \eqref{eq_2.1} we consider 
the Faddeev functions $G$, $\psi$, $h$ (see \cite{Faddeev1965}, \cite{Faddeev1974}, \cite{Henkin1987}, \cite{Novikov 1988}):
\begin{equation}\label{4.2}
	G(x,k) = e^{ikx} g(x,k), \ \ \ g(x,k) = - (2\pi)^{-3} \int\limits_{\mathbb{R}^3} \frac{e^{i\xi x} d\xi}{\xi^2 + 2k\xi},
\end{equation}
\begin{equation}\label{4.1}
	\psi(x,k) = e^{ikx} + \int\limits_{\mathbb{R}^3} G(x-y,k) v(y)\psi(y,k) dy, 
\end{equation}
where $x\in \mathbb{R}^3$, $k\in \mathbb{C}^3$,  $k^2=E$, $\mbox{Im}\, k\neq 0$,
\begin{equation}\label{4.4}
	h(k,l) = (2\pi)^{-3} \int\limits_{\mathbb{R}^3} e^{-ilx}v(x) \psi(x,k) dx,
\end{equation}
where 
\begin{equation}\label{4.5}
	 k,l \in \mathbb{C}^3, \ k^2=l^2=E,\ \mbox{Im}\, k = \mbox{Im}\, l \neq 0.
\end{equation}
One can consider (\ref{4.1}), (\ref{4.4}) assuming that 
\begin{equation}\label{4.6}
	\begin{aligned}
	v \text{ is a sufficiently regular function on } \mathbb{R}^3 \\
	\text{ with suffucient decay at infinity.}
	\end{aligned}
\end{equation}
For example, in connection with Problems 1.1 and 1.2, one can consider (\ref{4.1}), (\ref{4.4}) assuming that
\begin{equation}\label{4.7}
	v \in \mathbb{L}^{\infty}(B_{r_1}), \ \ \ v \equiv 0 \text{ on } \mathbb{R}^3\setminus B_{r_1}.
\end{equation}

We recall that (see \cite{Faddeev1965}, \cite{Faddeev1974}, \cite{Henkin1987}, \cite{Novikov 1988}): 
\begin{itemize}
\item The function $G$ satisfies the equation
\begin{equation}\label{4.8}
	(\Delta+E) G(x,k) = \delta(x), \ \ x\in\mathbb{R}^3, \ \ k \in \mathbb{C}^3\setminus \mathbb{R}^3, \ \ E=k^2;
\end{equation}

\item Formula (\ref{4.1}) at fixed $k$ is considered as an equation for 
\begin{equation}
	\psi = e^{ikx}\mu(x,k),
\end{equation}
where $\mu$ is sought in $\mathbb{L}^{\infty}(\mathbb{R}^3)$; 

\item As a corollary of (\ref{4.1}), (\ref{4.2}), (\ref{4.8}), 
$\psi$ satisfies (\ref{eq_2.1}) for $E=k^2$; 

\item The Faddeev functions $G$, $\psi$, $h$ are (non-analytic) continuation to the complex domain of functions    
 of the classical scattering theory for the Schr\"odinger equation (in particular, $h$ is a generalized "`scattering"' amplitude). 
\end{itemize} 
 
 In addition, $G$, $\psi$, $h$ in their zero energy restriction, that is for $E=k^2=0$, were considered for the first time in \cite{Beals1985}.
The Faddeev functions $G$, $\psi$, $h$ were, actually, rediscovered in \cite{Beals1985}.

Let 
\begin{equation}
\begin{aligned}
		\Sigma_E = \left\{ k\in \mathbb{C}^3: k^2 = k_1^2 + k_2^2 + k_3^2 = E\right\},\\
		\Theta_E = \left\{ k\in \Sigma_E,\  l\in\Sigma_E: \mbox{Im}\, k = \mbox{Im}\, l\right\},\\
		|k| = (|\mbox{Re}\,k|^2 +|\mbox{Im}\,k|^2)^{1/2}.
\end{aligned}		
\end{equation}
Let
\begin{equation}\label{2.5}
	\begin{aligned}
	\mathbb{W}^{m,q}(\mathbb{R}^3) = \{w:\  \partial^J w \in \mathbb{L}^q(\mathbb{R}^3),\  |J| \leq m \},\ \ m\in \mathbb{N}\cup 0,\ q\geq 1, \\
J \in (\mathbb{N}\cup 0)^3,\ |J| = \sum\limits_{i=1}\limits^{3}J_i,\ \partial^J v(x) 
= \frac{\partial^{|J|} v(x)}{\partial x_1^{J_1}\partial x_2^{J_2} \partial x_3^{J_3}},\\
	||w||_{m,q} = \max\limits_{|J|\leq m} ||\partial^J w||_{\mathbb{L}^q(\mathbb{R}^3)}.
	\end{aligned}
\end{equation}
Let the assumptions of Theorems 1.1 and 1.2 be fulfilled:
\begin{equation}\label{assumption}
		\begin{aligned}
		 (1-n) \in \mathbb{W}^{m,1}(\mathbb{R}^3) \text{ for some } m > 3, \\ 
		 \mbox{Im}\, n(x) \geq 0, \ \ \ x\in \mathbb{R}^3,\\
		 \mbox{supp}\, (1 - n) \subset B_{r_1},  \\
		 \|1 - n\|_{m,1} \leq C_n. 
	\end{aligned}
\end{equation}
Let
		\begin{equation}
			v=\omega^2(1-n),\ \ \ N = \omega^2 C_n, \ \ \ E = \omega^2.
		\end{equation}
Then we have that:
\begin{equation}\label{mu_1}
	\mu(x,k) \rightarrow 1 \ \  \text{ as } \ \ |k|\rightarrow \infty 
\end{equation}
and, for any $\sigma>1$,
\begin{equation}\label{mu_2}
	|\mu(x,k)|  \leq \sigma \ \  \text{ for } \ \ |k| \geq \lambda_1(N, m,\sigma, r_1),
\end{equation}
where $x\in \mathbb{R}^3$, $k \in \Sigma_E$;
\begin{equation}\label{lim_1}
			\hat{v}(p) = \lim\limits_
		{\scriptsize
			\begin{array}{c}
			(k,l)\in \Theta_E,\, k-l=p\\
			|\mbox{Im}\,k|=|\mbox{Im}\,l|\rightarrow \infty
			\end{array}
		} h(k,l)\ \ \ 	 \text{ for any } p\in \mathbb{R}^3,
	\end{equation}
\begin{equation}\label{lim_2}
	\begin{aligned}
		|\hat{v}(p) - h(k,l)|\leq \frac{c_1(m,r_1)N^2}{(E+\rho^2)^{1/2}}  	\ \ \text{ for }  (k,l) \in \Theta_E, \ \  p = k-l,\\
		|\mbox{Im}\,k| = |\mbox{Im} \,l| = \rho, \ \ \ E+\rho^2 \geq \lambda_2(N,m,r_1), \\ p^2 \leq 4(E+\rho^2),
	\end{aligned}
\end{equation}
		where 
		\begin{equation}
			\hat{v}(p) = (2\pi)^{-3}\int\limits_{\mathbb{R}^3} e^{ipx} v(x)dx, \ \ p\in \mathbb{R}^3.
		\end{equation}
Results of the type (\ref{mu_1}), (\ref{mu_2}) go back to \cite{Beals1985}.  
For more information concerning (\ref{mu_2}) see estimate (4.11) of \cite{IN2012}.
Results of the type (\ref{lim_1}), (\ref{lim_2}) 
(with less precise right-hand side in (\ref{lim_2})) go back to \cite{Henkin1987}. 
Estimate (\ref{lim_2}) follows, for example, from 
formulas (\ref{4.1}), (\ref{4.4}) and 
the estimate
\begin{equation}\label{4.15.new}
	\begin{aligned}
	\| \Lambda^{-s} g(k) \Lambda^{-s}\|_{\mathbb{L}^2(\mathbb{R}^d)\rightarrow \mathbb{L}^2(\mathbb{R}^d)} = O(|k|^{-1}) \\ \text{ as } \ |k|\rightarrow \infty,\ \ \
		k\in \mathbb{C}^3\setminus \mathbb{R}^3,
	\end{aligned}
\end{equation}
for $s>1/2$, where $g(k)$ denotes the integral operator with the Schwartz kernel $g(x-y,k)$ and $\Lambda$ denotes the multiplication operator by the function $(1+|x|^2)^{1/2}$. Estimate (\ref{4.15.new}) was formulated, first, in \cite{LN1987}.
This estimate generilizes, in particular, some related estimate of \cite{SU1987} for $k^2=E=0$.
 Concerning proof of (\ref{4.15.new}), see \cite{Weder1991}. 

In addition, we have that:
\begin{equation}\label{delta_h0}
	\begin{aligned}
	h_2(k,l) - h_1(k,l) = (2\pi)^{-3} \int\limits_{\mathbb{R}^3} \psi_1(x,-l) (v_2(x) - v_1(x)) \psi_2(x,k) dx 	\\ 	 
	\text{ for }  (k,l) \in \Theta_E,\ |\mbox{Im}\,k| = |\mbox{Im} \,l| \neq 0,\\ \text{ and $v_1$, $v_2$ satisfying (\ref{4.6}),}
	\end{aligned}
\end{equation}
and, under the assumptions of Theorems \ref{Theorem_1.1} and \ref{Theorem_1.2},
\begin{equation}\label{delta_v}
	\begin{aligned}
		|\hat{v}_1(p) - \hat{v}_2(p) - h_1(k,l) + h_2(k,l)|
		\leq \frac{
		c_2(m,r_1)N
		\|v_1-v_2\|_{\mathbb{L}^{\infty}(B_{r_1})}
		}
		{(E+\rho^2)^{1/2}}  	\\ \text{ for }  (k,l) \in \Theta_E, \ \  p = k-l,\ \
		|\mbox{Im}\,k| = |\mbox{Im} \,l| = \rho,\\ E+\rho^2 \geq \lambda_3(N,m,r_1), \ \ p^2 \leq 4(E+\rho^2),
	\end{aligned}
\end{equation}
where $h_j$, $\psi_j$ denote $h$ and $\psi$ of (\ref{4.4}) and (\ref{4.1}) for 
$v_j = \omega^2(1-n_j)$,  $j=1,2$, $N=\omega^2 C_n$, $E=\omega^2$.

Formula (\ref{delta_h0}) was given in \cite{Novikov1996}, \cite{Novikov2005}.
Estimate \eqref{delta_v} was given e.g. in \cite{IN2012++}.

\section{Proofs  of Theorem \ref{Theorem_1.1} and Theorem 1.2}
{\it3.1. Preliminaries.}
In this section we always assume for simplicity that $r_1=1$.

We consider the operators $\hat{\mbox{S}}_{j}$, $j=1,2$, defined as follows
\begin{equation}
	(\hat{\mbox{S}}_{j} \phi) (x) = \int\limits_{\partial B_r} G^+_{j} (x,y,\omega) \phi(y) dy, \ \ \ x\in \partial B_r, \ j=1,2.
\end{equation}
Note that
\begin{equation}\label{eq4.2}
	\|\hat{\mbox{S}}_{1} - \hat{\mbox{S}}_{2} \|_{\mathbb{L}^2(\partial B_r)} \leq 
	\|G^+_{1} - G^+_{2}  \|_{\mathbb{L}^2(\partial B_r) \times \mathbb{L}^2(\partial B_r)}.
\end{equation}
To prove Theorems 1.1 and 1.2 we use, in particular, the following lemmas (see Lemma 3.2 
and proof of Theorem 1.2 of \cite{HH2001}):
\begin{Lemma}\label{Lemma_3.1}
	 Assume $r_1=1<r<r_2$. Moreover,  $n_1$, $n_2$ are refractive indices
	 with {\rm $\mbox{supp}\, (1-n_1), \mbox{supp}\, (1-n_2)\subset B_{1}$}. Then, there exists  
	a postive constant $c_3$ (depending only on $\omega$, $r$, $r_2$) such that for all solutions 
	$\psi_1\in C^2(B_{r_2}) \cap \mathbb{L}^2(B_{r_2})$ to $\Delta\psi + \omega^2n_1\psi = 0$ in $B_{r_2}$
	and  all solutions 
	$\psi_2\in C^2(B_{r_2}) \cap \mathbb{L}^2(B_{r_2})$ to $\Delta\psi + \omega^2n_2\psi = 0$ in $B_{r_2}$
	the following estimate holds:
	\begin{equation}\label{eq4.3}
		\begin{aligned}
		\left|\int\limits_{B_1} (n_1-n_2)\psi_1\psi_2 dx\right| \leq c_3
		\|\hat{\mbox{S}_{1}} - \hat{\mbox{S}_{2}} \|_{\mathbb{L}^2(\partial B_r)} \|\psi_1\|_{\mathbb{L}^2(B_{r_2})} 
		\|\psi_2\|_{\mathbb{L}^2(B_{r_2})}. 
		\end{aligned}
	\end{equation}
\end{Lemma}
Note that estimate \eqref{eq4.3} is derived in \cite{HH2001} using an Alessandrini type identity,
where instead of the Dirichlet-to-Neumann maps the operators $\hat{\mbox{S}}_{1}$, $\hat{\mbox{S}}_{2}$ are used, 
see \cite{Alessandrini1988}, \cite{HH2001}.

\begin{Lemma}\label{Lemma_3.2}
	 Let $r>r_1=1$, $\omega>0$, $C_n>0$, $\mu>3/2$ and $0<\theta<1$. Let $n_1,n_2$ be  refractive indices such that
	 $\|(1-n_j)\|_{\mathbb{H}^{\mu}(\mathbb{R}^3)} \leq C_n$, {\rm $\mbox{supp} (1-n_j)\subset B_{1}$}, $j=1,2$, where $\mathbb{H}^{\mu} = \mathbb{W}^{\mu,2}$. 
	 Then there exist positive constants $T$ and $\eta$ such that
	\begin{equation}\label{GGff}
		\begin{aligned}
		\|G_1^+-G_2^+\|_{\mathbb{L}^2(\partial B_{2r} \times\partial B_{2r} )}^2 
		\leq \eta^2 \exp\left(-\left( -\ln \frac{\|f_1-f_2\|_{\mathbb{L}^2({\mathcal M}_\omega)}}{T\eta}\right)^\theta \right)
		\end{aligned}
	\end{equation}
	for sufficiently small $\|f_1-f_2\|_{\mathbb{L}^2({\mathcal M}_\omega)}$, where $G_j^+$, $f_j$ are near and far field
	scattering data for $n_j$, $j=1,2$, at fixed frequency $\omega$.
\end{Lemma}

\noindent
{\it 3.2. Proof of Theorem 1.1.} Let 
\begin{equation}\label{6.3}
\begin{aligned}
	\mathbb{L}^{\infty}_{\mu}(\mathbb{R}^3) = \{u\in \mathbb{L}^{\infty}(\mathbb{R}^3): \|u\|_\mu<+\infty\},\\
	\|u\|_\mu = \mbox{ess}\,\sup\limits_{p\in\mathbb{R}^3} (1+|p|)^{\mu}|u(p)|, \ \ \ \mu>0.
\end{aligned}
\end{equation}
Note that
\begin{equation}\label{6.4}
	\begin{aligned}
	w \in \mathbb{W}^{m,1}(\mathbb{R}^3) \Longrightarrow \hat{w} \in 	\mathbb{L}^{\infty}_{\mu}(\mathbb{R}^3)\cap { C}(\mathbb{R}^3),\\
		\|\hat{w}\|_\mu \leq c_4(m) \|w\|_{m,1} \ \ \ \text{ for } \ \ \mu = m, 
	\end{aligned}
\end{equation}
where $\mathbb{W}^{m,1}$,  $\mathbb{L}^{\infty}_{\mu}$ are the spaces of (\ref{2.5}), (\ref{6.3}),
\begin{equation}
	\hat{w}(p) = (2\pi)^{-3}\int\limits_{\mathbb{R}^3} e^{ipx}w(x) dx, \ \ \ p\in \mathbb{R}^3.
\end{equation}

Let
\begin{equation}\label{eq3.7}
 N = \omega^2 C_n, \ \ \ E = \omega^2,\ \ \ 	v_j=\omega^2(1-n_j),\ j=1,2.
\end{equation}
Using the inverse Fourier transform formula
\begin{equation}
	w(x) = \int\limits_{\mathbb{R}^3} e^{-ipx}\hat{w}(p) dp, \ \ \ x\in \mathbb{R}^3,
\end{equation}
we have that
\begin{equation}\label{5.4}
	\begin{aligned}
		\|v_1 - v_2\|_{\mathbb{L}^{\infty}(D)} \leq 
		\sup\limits_{x\in \bar{B}_1} 
		\left|\int\limits_{\mathbb{R}^3} e^{-ipx}\left(\hat{v}_2(p)
		 - \hat{v}_1(p)\right) dp\right| \leq\\
		\leq I_1(\kappa) + I_2(\kappa) \ \ \ \text{ for any }  \ \kappa>0,
	\end{aligned}
\end{equation}
where
\begin{equation}\label{6.7}
	\begin{aligned}
		I_1(\kappa) = \int\limits_{|p|\leq \kappa} |\hat{v}_2(p) - \hat{v}_1(p)| dp, \\
		I_2(\kappa) = \int\limits_{|p|\geq \kappa} |\hat{v}_2(p) - \hat{v}_1(p)| dp.
	\end{aligned}
\end{equation}
Using (\ref{6.4}), we obtain that
\begin{equation}\label{6.8}
	|\hat{v}_2(p) - \hat{v}_1(p)| \leq 2c_4(m) N (1+|p|)^{-m}, \ \ \ p\in \mathbb{R}^3.
\end{equation}
Using (\ref{6.7}), (\ref{6.8}), we find that, for any $\kappa>0$,
\begin{equation}\label{5.16}
	\begin{aligned}
		I_2(\kappa) \leq 8\pi c_4(m) N  \int\limits_{\kappa}\limits^{+\infty} \frac{dt}{t^{m-2}} \leq 
		\frac{8\pi c_4(m)N}{m-3}\frac{1}{\kappa^{m-3}}.
	\end{aligned}
\end{equation}

Due to (\ref{delta_v}), we have that
\begin{equation}\label{6.9}
\begin{aligned}
	|\hat{v}_2(p) - \hat{v}_1(p)| \leq |h_2(k,l) - h_1(k,l)| + 
	\frac{c_2(m)N\|v_1-v_2\|_{\mathbb{L}^{\infty}(B_1)}}{(E+\rho^2)^{1/2}},\\
	\text{ for }  (k,l) \in \Theta_E, \ \  p = k-l,\ \
		|\mbox{Im}\,k| = |\mbox{Im} \,l| = \rho,\\ E+\rho^2 \geq \lambda_3(N,m), \ \ p^2 \leq 4(E+\rho^2).
\end{aligned}
\end{equation}
Let
\begin{equation}
	\begin{aligned}
	r_2 \text{ be some fixed constant such that } r_2>r, \\  
	\delta = ||G^+_{1} - G^+_{2}||_{\mathbb{L}^{2}(\partial B_r \times \partial B_r)},\\
	c_5 = (2\pi)^{-3}\int\limits_{ B_{r_2}} dx.
	\end{aligned}
\end{equation}
Combining  (\ref{delta_h0}), \eqref{eq4.2}, \eqref{eq4.3} and \eqref{eq3.7}, we get that
\begin{equation}\label{6.11}
	\begin{aligned}
	|h_2(k,l) - h_1(k,l)| &\leq \\ \leq  c_3 c_5 \omega^2
	&\|\psi_1(\cdot,-l)\|_{\mathbb{L}^{\infty}(B_{r_2})}\, \delta \,
	\|\psi_2(\cdot,k)\|_{\mathbb{L}^{\infty}(B_{r_2})},\\
	&(k,l)\in \Theta_E, \ |\mbox{Im}\,k| = |\mbox{Im}\,l| \neq 0.
	\end{aligned}
\end{equation}
Using (\ref{mu_2}), we find that
\begin{equation}\label{6.12}
\begin{aligned}
	\|\psi_j(\cdot,k)\|_{\mathbb{L}^{\infty}(B_{r_2})} \leq
	 \sigma \,\exp\bigg(|\mbox{Im}\,k|r_2\bigg), \ \
	j=1,2, \\ k \in \Sigma_E, \ |k|\geq \lambda_1(N,m,\sigma).
\end{aligned}
\end{equation}
Here and bellow in this section the constant $\sigma$ is the same that in (\ref{mu_2}).
 
Combining (\ref{6.11}) and (\ref{6.12}), we obtain that
\begin{equation}\label{6.13}
	\begin{aligned}
		|h_2(k,l) - h_1(k,l)| \leq c_3c_5 \omega^2\sigma^2  e^{2\rho r_2} \delta,\\ \text{ for }  
		(k,l)\in\Theta_E, \ \ \rho = |\mbox{Im}\,k|= |\mbox{Im}\,l|, \\ E+\rho^2\geq \lambda_1^2(N,m,\sigma).
	\end{aligned}
\end{equation}
Using (\ref{6.9}), (\ref{6.13}), we get that
\begin{equation}\label{6.14}
	\begin{aligned}
		|\hat{v}_2(p) - \hat{v}_1(p)| \leq c_3 c_5 \omega^2\sigma^2   e^{2\rho r_2} \delta 
		+ \frac{c_2(m)N\|v_1-v_2\|_{\mathbb{L}^{\infty}(B_1)}}{(E+\rho^2)^{1/2}},\\
		p\in\mathbb{R}^3,\  p^2 \leq 4(E+\rho^2),\  E+\rho^2 \geq \max\{ \lambda_1^2 , \lambda_3\}.
	\end{aligned}
\end{equation}
Let 
\begin{equation}
	\varepsilon = \left(\frac{3}{8\pi c_2(m)N}\right)^{1/3}
\end{equation}
and $\lambda_4(N,m,\sigma)>0$ be such that
\begin{equation}\label{5.14new.}	
	E+\rho^2 \geq \lambda_4(N,m,\sigma) \Longrightarrow 
	\left\{
	\begin{aligned}
	&E+\rho^2 \geq \lambda_1^2(N,m,\sigma), \\ &E+\rho^2 \geq \lambda_3(N,m), \\ &
	\left(\varepsilon(E+\rho^2)^{\frac{1}{6}}\right)^2 \leq 4(E+\rho^2).
	\end{aligned}\right.
\end{equation} 
Using (\ref{6.7}), (\ref{6.14}), we get that
\begin{equation}\label{5.15}
\begin{aligned}
	I_1(\kappa) \leq \frac{4}{3} \pi\kappa^3 \Big( c_3 c_5 \omega^2\sigma^2   e^{2\rho r_2} \delta 
		+ \frac{c_2(m)N\|v_1-v_2\|_{\mathbb{L}^{\infty}(B_1)}}{(E+\rho^2)^{1/2}}\Big),\\  
		\kappa>0, \ \kappa^2\leq 4(E+\rho^2), \\ 
	E+\rho^2 \geq \lambda_4(N,m,\sigma).
\end{aligned}
\end{equation}
 Combining (\ref{5.4}), (\ref{5.16}), (\ref{5.15})  for $\kappa = \varepsilon(E+\rho^2)^{\frac{1}{6}}$ and (\ref{5.14new.}), we get that
 \begin{equation}\label{5.18}
 \begin{aligned}
 	\|v_1 - v_2\|_{\mathbb{L}^{\infty}(B_1)} \leq c_6(N,m,\omega,\sigma) \sqrt{E+\rho^2}\, e^{2\rho r_2} \delta+ \\+ 
 	c_7(N,m)(E+\rho^2)^{-\frac{m-3}{6}} + \frac{1}{2}\|v_1 - v_2\|_{\mathbb{L}^{\infty}(B_1)}, \\
 	E+\rho^2 \geq \lambda_4(N,m,\sigma).
 	\end{aligned}
 \end{equation}
 
Let $\tau \in (0,1)$ and 
\begin{equation}\label{5.19}
	\beta = \frac{1-\tau}{2r_2}, \ \ \ \rho = \beta\ln\left(3+ \delta^{-1}\right),
\end{equation} 
where $\delta$ is so small that $E+\rho^2 \geq \lambda_4(N,m,\sigma)$. Then due to (\ref{5.18}), we have that
\begin{equation}\label{5.20}
	\begin{aligned}
		\frac{1}{2}\|v_1 - v_2\|_{\mathbb{L}^{\infty}(D)} &\leq \\ \leq 
		c_6(N,m,\omega,\sigma)  &\left(E+ \left(\beta\ln\left(3+ \delta^{-1}\right)\right)^2\right)^{1/2}
		\left(3+ \delta^{-1}\right)^{2\beta r_2} \delta + \\
		+  c_7(N,m) &\left(E+ \left(\beta\ln\left(3+ \delta^{-1}\right)\right)^2\right)^{-\frac{m-3}{6}} =\\
		= c_6(N,m,\omega,\sigma)  &\left(E+ \left(\beta\ln\left(3+ \delta^{-1}\right)\right)^2\right)^{1/2}
		\left(1+ 3\delta \right)^{1-\tau} \delta^{\tau} 
		 + \\
		&+ c_7(N,m) \left(E+ \left(\beta\ln\left(3+ \delta^{-1}\right)\right)^2\right)^{-\frac{m-3}{6}},
	\end{aligned}
\end{equation}
where  $\tau, \beta$ and $\delta$ are the same as in (\ref{5.19}). 

Using (\ref{5.20}), we obtain that
\begin{equation}\label{5.21}
	\|v_1 - v_2\|_{\mathbb{L}^{\infty}(B_1)} \leq c_{8}(N,m,\omega,\sigma,\tau) \left(\ln\left(3+\delta^{-1}\right)\right)^{-\frac{m-3}{3}}
\end{equation}
for $\delta = ||G^+_{1} - G^+_{2}||_{\mathbb{L}^{2}(\partial B_r \times \partial B_r)}\leq \delta_1(N,m,\omega,\sigma,\tau)$, where
$\delta_1$ is a sufficiently
small positive constant. Estimate (\ref{5.21}) in the general case (with modified $c_{8}$) follows from
 (\ref{5.21})  for $\delta \leq \delta_1(N,m,\omega,\sigma,\tau)$ 
and the property that 
\begin{equation}\label{new23}
\|v_j\|_{\mathbb{L}^{\infty}(B_1)} \leq c_{9}(m)N, \ \ j=1,2.
\end{equation}
Taking into account \eqref{eq3.7}, we obtain \eqref{eq_t1}.

\noindent
{\it 3.2. Proof of Theorem 1.2.}
According to the Sobolev embedding theorem, we have that
\begin{equation}\label{eq_emb}
	\mathbb{W}^{m,1}(\mathbb{R}^3) \subset \mathbb{H}^{m-3/2}(\mathbb{R}^3),
\end{equation}
where $\mathbb{H}^{\mu} = \mathbb{W}^{\mu,2}$. 

Combining  \eqref{eq_c},  \eqref{eq_t1}, \eqref{GGff} 
with $\theta$ satisfying $\theta \frac{m-3}{3} = \frac{m-3}{3} - \epsilon$, 
and \eqref{eq_emb}, we obtain \eqref{eq_t2} for sufficiently small $\|f_1-f_2\|_{\mathbb{L}^2({\mathcal M}_\omega)}$ (analogously with the proof of Theorem 1.2 of \cite{HH2001}).
Using also \eqref{new23} and \eqref{eq3.7}, we get estimate \eqref{eq_t2} in the general case.


\section*{Acknowledgements}
This work was partially supported by FCP Kadry No. 14.A18.21.0866, Moscow Institute of Physics and Technology.

\noindent
{ {\bf M.I. Isaev}\\
Centre de Math\'ematiques Appliqu\'ees, Ecole Polytechnique,

91128 Palaiseau, France\\
Moscow Institute of Physics and Technology,

141700 Dolgoprudny, Russia\\
e-mail: \tt{isaev.m.i@gmail.com}}\\

\noindent
{ {\bf R.G. Novikov}\\
Centre de Math\'ematiques Appliqu\'ees, Ecole Polytechnique,

91128 Palaiseau, France\\
Institute of Earthquake Prediction Theory and Math. Geophysics RAS,

117997 Moscow, Russia\\ 
e-mail: \tt{novikov@cmap.polytechnique.fr}}

\end{document}